\documentclass[submission,copyright,creativecommons]{eptcs}
\providecommand{\event}{ThEdu 2018 Postproceedings} 
\providecommand{\thisvolume}[1]{this volume of EPTCS, Open Publishing Association}
\usepackage{underscore}           
\usepackage{paralist}       
\usepackage{graphicx}       
\usepackage{amsmath}        
\usepackage{hyperref}  
\usepackage{breakurl}  

\title{Technologies for \\``Complete, Transparent \& Interactive Models of Math'' \\ in Education}
\author{
Walther Neuper
\institute{University of Technology\\
Graz, Austria}
\email{wneuper@ist.tugraz.at}
}
\def\titlerunning{Technologies for Models of Mathematics}
\def\authorrunning{W.Neuper}

\begin{document}
\maketitle

\begin{abstract}
A new generation of educational mathematics software is being shaped in ThEdu~\footnote{\url{https://www.uc.pt/en/congressos/thedu/}} and other academic communities on the side of computer mathematics. Respective concepts and technologies have been clarified to an extent, which calls for cooperation with educational sciences in order to optimise the new generation's impact on educational practice.
The paper addresses educational scientists who want to examine specific software features and estimate respective effects in STEM education at universities and subsequently at high-school.

The key features are characterised as a ``complete, transparent and interactive model of mathematics'', which offers interactive experience in all relevant aspects in doing mathematics. Interaction uses several layers of formal languages: the language of terms, of specifications, of proofs and of program language, which are connected by Lucas-Interpretation providing ``next-step-guidance'' as well as providing prover power to check user input.

So this paper is structured from the point of view of computer mathematics and thus cannot give a serious description of effects on educational practice --- this is up to collaboration with educational science; such collaboration is prepared by a series of questions, some of which are biased towards software usability (and mainly to be solved by computer mathematicians) and some of which are biased towards genuine research in educational sciences.
\end{abstract}

\section{Introduction}\label{sec:intro}
The quality of formula-oriented\footnote{Software modelling geometry and other non-formal mathematics is excluded from consideration in this paper.} software available for mathematics education is far beyond the state of the art in computer mathematics, in particular of the state in (computer) theorem proving (abbreviated ``TP''). The state in respective technologies is being presented and discussed mainly in conferences and publications which do not reach educators --- so there is the issue to support interdisciplinary communication between disciplines, who need to cooperate in development of educational software towards the state of the art.

This paper intends to reconfirm communication between TP experts and educational sciences. The latter are facing educational software for mathematics like sand on the sea and they hardly believe, that there can be software completely different from existing products; and educators are also aware, that much software is part of the problem in mathematics education and not a solution. 

The paper originates from thorough discussions under different aspects at two conferences~\footnote{\url{http://www.ist.tugraz.at/projects/isac/publ/mechanise-expl-prensent.pdf}}
\footnote{\url{http://www.ist.tugraz.at/projects/isac/publ/justify-sys-explain-present.pdf}}, which aim at interdisciplinarity and inclusion of education themselves. Particularly influential were discussions at a recent workshop on semiotics~\footnote{\url{http://wwwu.uni-klu.ac.at/kadunz/semiotik/news.html}} which allowed to essentially refine and complete an extended abstract \cite{wn:cme-ei-18}.
 
\medskip
Following the point of view in semiotics, mathematics is regarded as a language; thus here technologies of computer mathematics are structured alongside languages implemented here. The paper will show that specific kinds of languages combined in a novel conceptual view as different language layers can implement ``complete, transparent and interactive'' models of mathematics. ``Model of mathematics'' indicates a technical approach to educational software: The first step does \emph{not} consider educational theory, \emph{not} tutoring, \emph{not} instructional design and the like --- just \emph{software reflecting the structure of mathematics}~\footnote{During the last century mathematics has been formalised to an extent that allows to implement major aspects of doing mathematics in computer software --- to an extent far beyond other scientific disciplines, all of which are more immediately related to physical matter than mathematics (where physics cannot be implemented, only respective symbols).} in a model featuring ``completeness, transparency and interactivity'', features to be described in detail below.

The paper regards ``mathematics'' as what is taught as ``mathematics for \dots'' in science, technology and engineering studies~\footnote{The abbreviation STEM includes mathematics, however, studies in pure mathematics are excluded from this paper.}. This restriction, in particular, excludes explicit teaching to prove, because engineers like to rely on the correctness of what they adopt from mathematicians. And this restriction is conforming to common practice in high-school mathematics --- however, the paper will also show, that TP involves rigorous formal proof and educational software should be based on that behind the scenes; and thus effects can be expected on students working on such technology, even if they do \emph{not} explicitly use that software for proving.

\medskip
Technologies of computer mathematics as addressed in this paper are rapidly evolving. In order to achieve a reliable source for educators, who want to be informed of these technologies and their relevance for mathematics education, the paper separates technical facts (already implemented in prototypes) and prospects (already in a concrete stage of planning for different software products) from foreseeable effects in education. Examples will be mostly taken from a prototype development~\cite{wn:proto-sys}~\footnote{\url{http://www.ist.tugraz.at/isac/History}} where the author is involved and early field tests were successful \cite{imst-htl06,imst-htl07,imst-hpts08}.

\paragraph{The paper is structured} as follows: The subsequent four sections each describe a language layer, each separated into technical facts / prospects and foreseeable software features / effects on education together with questions, where \textbf{?U} marks questions mainly addressing software usability and \textbf{?R} marks questions mainly addressing genuine educational research\label{pt-U-R}: 
\S\ref{sec:term-lang} the language of mathematical terms --- and transparency as regards definitions, types, etc, 
 \S\ref{sec:spec-lang} the language layer of specifications --- and respective means to ``divide and conquer'' problems, \S\ref{sec:proof-lang} the ``proof language'' --- adapted to habits and needs of engineers and \S\ref{sec:prog-lang} the programming language required to implement algorithms solving classes of problems. Another section, \S\ref{sec:next-step}, describes a technology crossing over all four language layers and providing ``next-step-guidance'' to students who get stuck during problem solving --- a feature indispensable for ``systems that explain themselves'', an aim addressed in \S\ref{sec:summary}, the summary and conclusion.

\section{Term Language --- Transparent to Definitions}\label{sec:term-lang}
The language of terms comprises propositions, predicates and algebraic expressions; it is the bottom of other languages in formal mathematics. 

\subsection{Technical Facts and Prospects}\label{ssec:term-facts}
In this section the proof assistant Isabelle \cite{Isa-tutor08} will serve as an example in technical details, in the implementation already existing as well as in planning future features.

\paragraph{Terms are typed:} Here the notion ``type'' refers to simple typed lambda calculus \cite{pl:barend85}, which is completely different from what is known as ``type'' in programming languages --- a frequent source of misunderstanding in communication with educators and educational scientists. In TP a type, for instance the integer numbers, is not defined as a certain bit-pattern like in programming, but as a set of objects gathered by a formal definition --- in Fig.\ref{fig:isa-transparent} below the definition \texttt{intrel\_iff} shows exactly what mathematics students are taught in introductory courses: integers defined as equivalence classes on pairs of natural numbers $(x,y)\approx(u,v)\Longleftrightarrow x+v = u+y$.

\paragraph{Definitions are reachable by mouse-click} (as well as types mentioned above): What is an indispensable feature of a programmers' IDEs and a prerequisite to master the complexity of contemporary software, that is now also available in most advanced IDEs for proof assistants; Fig.\ref{fig:isa-transparent} shows a compilation of screenshots from Isabelle/jEdit~\cite{DBLP:journals/corr/Wenzel14}:
\begin{figure} [htb]
  \centering
  \includegraphics[width=100mm]
    {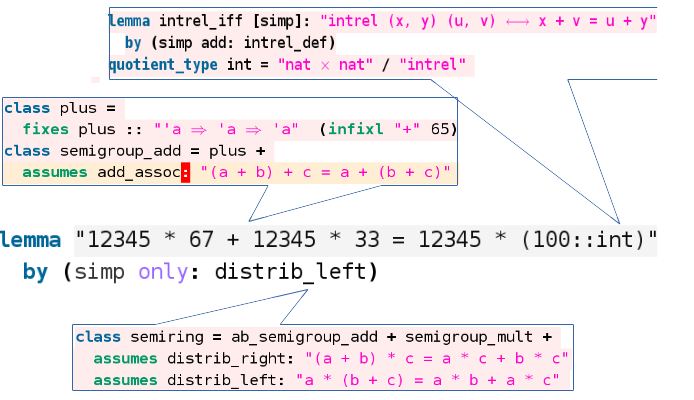}
  \caption{Isabelle shows underlying knowledge transparently.}
  \label{fig:isa-transparent}
\end{figure}
Given a simple calculation as a \texttt{lemma} (as a simple example for more comprehensive theorems) a mouse-click on an element immediately leads to the respective definition; a click on \texttt{int} leads to the definition mentioned above, a click on the \texttt{+} leads to the definition close to relevant algebraic structures like the additive semi-group \texttt{semigroup\_add}, a click on the theorem \texttt{distrib\_left} to be used by the \texttt{simp}lifier leads to the respective formula in the context of a \texttt{semiring}, where the integers in the calculation at hand are an instantiation of (which is formally proved also there). Other semantic information is available as well, for instance distinction of constants, of free and bound variables by colours.

Expert-users of Isabelle know, that not (yet) all definitions are reachable by mouse-click: In contrary to programs in software development, not all knowledge of Isabelle can be loaded into the computer's memory. The current Isabelle version requires 16GB RAM, but that is far too small to load the standard distribution's knowledge~\footnote{\url{https://isabelle.in.tum.de/website-Isabelle2018/dist/library/HOL/index.html}},
not to speak of Isabelle's archive of formal proofs~\footnote{\url{https://www.isa-afp.org/}}. So, there are plans to store all linked data created during theory processing in a database, which can be referenced via networks --- this is an answer anticipating educators' questions, why an Isabelle-based educational software, whether on a handheld device or not,  can only be available as a cloud service and can\emph{not} run off-line.

\paragraph{Presentation of terms} and respective editing is according to requirements of software engineers (in verification tasks) and of mathematicians (mechanising mathematical definitions and proofs). The Isabelle/jEdit IDE is the most advanced among the leading proof assistants; it presents terms nicely with exponents, subscripts and many special characters like $\sum,\int,\Rightarrow,\forall,\exists,\subseteq,\in,\otimes,\sqcap$ etc.~--- and this presentation is open for adding arbitrary operators required in future extensions of knowledge.

However, this presentation is still fixed to lines and bypasses 2-dimensional presentation as known, for instance, from computer algebra systems; $\frac{x}{2}$ is presented as $x/2$ which makes more complicated fractions unreadable --- which is unacceptable for engineers and for education. This is an issue to be tackled from the side of development projects for educational software \emph{alone}, there appears no interest from the side of proof assistants --- their users are satisfied with the present situation. 

The technical side of the issue is to develop formula editors, which combine the 2-dimensional representation of terms as known from algebra systems with the ability to link elements of terms with respective semantic information as known from Isabelle. A feasibility study~\cite{mmahringer} indicates, that development efforts to combine the two features are manageable.

\subsection{Foreseeable Effects and Questions}\label{ssec:term-effects}
Several Effects of the above-mentioned technical facts and plans can be foreseen. The three most striking effects are listed below; their impact on education, as said already, can\emph{not} be anticipated easily; so it is addressed by questions, which are to be answered in collaboration between computer mathematicians and by educational scientists and some of which seem new within educational sciences~\footnote{\textbf{?U} marks questions mainly addressing software usability and \textbf{?R} marks questions mainly addressing educational research as already announced in \S\ref{sec:intro} on p.\pageref{pt-U-R}}.

\paragraph{All relevant mathematical knowledge is accessible by a student} who uses a TP-based system, and the access is as easy and quick as a software-engineer's access to program code. In TP such access concerns formal definitions, for instance that for \texttt{int} inf Fig.\ref{fig:isa-transparent}. Apparently formal definitions need to be combined with explanations (which are not present in proof assistants, because respective users already know how to read terms) for undergraduates and pupils.
\begin{itemize}
\item[\textbf{?U:}] \textit{How guard students against being overtaxed,} if all knowledge underlying a concrete situation in mathematical problem solving is available any time (see next question) by mouse-click? Subsequent questions: Which kind of informal (multimedia-based!) explanations have to be attached to which definitions at which levels of knowledge? How filter off overwhelming technicalities in mechanised mathematics?

\item[\textbf{?U:}] \textit{How lure students into interactive experiments with definitions?} The strength of terms implemented in TP-based systems is interactivity. There are promising developments to support students in interactive experiments with formal definitions in software verification \cite{DBLP:conf/mkm/Schreiner18} --- so the question is: How transfer such technologies to knowledge mechanised in proof assistants such that experimenting with definitions becomes instructive?
\end{itemize}

\paragraph{Knowledge is accessible any time} from any term: that means, all knowledge is accessible, as soon as a student has started to solve a problem by use of a term --- this term is the student's first entry to knowledge, and upon having reached the solution~\footnote{\S\ref{sec:next-step} will show, that constructing a solution of a problem is, in principle, possible in a completely passive way just by clicking.} she or he must have passed all knowledge required to solve the problem (in general without having noticed all that knowledge, which is waiting for the student's request).

\begin{itemize}
\item[\textbf{?R:}] \textit{Reconsider relations between competences and factual knowledge.} The presently prevailing theory of fundamental competences (see, for instance \cite{PISA-math-2015}) originates from outside mathematics and is increasingly debated as reason for decreasing mathematical knowledge of freshmen at universities (see \cite{stellungn-dmv-gdm-mnu} for Germany). Presence of mechanised knowledge, which is easier accessible than ever before, sharpens the question: What should students learn, solving problems (and again: what does that mean?), memorising definitions, clicking around in some software?

A few decades ago integration, for instance, was a task for demanded experts using Bronstein~\cite{bronstein-integral} and other tables, today everyone uses computer algebra for such tasks. So, what can be dropped in presence of TP-based educational systems~\cite{bb:box}?
\end{itemize}

\paragraph{A new kind of formula editor} \emph{must} be developed due to the challenging situation caused by TP's IDEs inappropriate for mathematics education. But meeting the challenge promises benefits if combined with features presently not available in computer algebra or in web-presentations. Some of them might be surprising, for instance: There are well-known methods to teach applying algebra rules~\cite{malle93:algebra}. However, these methods remain unused, because boring in practice without support. Now TP technology can provide that support (with typed terms even more accurately than computer algebra without). \cite{mmahringer} gives the example below, how a TP-based editor could interactively support a student applying the chain rule somewhere within a calculation (which will be introduced in \S\ref{sec:proof-lang} in detail, \texttt{04..08} are line numbers not belonging to the calculation, $\bigotimes$ represents the cursor requesting input, colours might further distinguish the boxes):

\medskip
{\footnotesize
\texttt{04} \dots\\

\texttt{05} 
  $\;\;\;\;\;\frac{d}{d\,x} x + \frac{d}{d\,x} \boxed{\sin(\boxed{x^2})}$\\

\texttt{06} \hspace{40mm} 
  ${\it Rewrite}\;(\frac{d}{d\,{\it bdv}} \boxed{\sin(\boxed{u})} = 
    \boxed{\cos(u)} * \frac{d}{d\,{\it bdv}} \boxed{u})$\\

\texttt{07} 
  $\;\;\;\;\;\frac{d}{d\,x} x + 
    \boxed{\cos(x^2)} * \frac{d}{d\,x} \boxed{\bigotimes}$\\

\texttt{08} \dots\\
} 

\noindent
Given such an editor these questions arise:
\begin{itemize}
\item[\textbf{?U:}] \textit{How exhibit students to these boxes without them getting bored?} In more technical detail: How design a user model, which remembers what abilities a student already has achieved? How balance student's passive clicking (probably in steps \texttt{05} and \texttt{06} above) and active editing (probably a whole differentiation next time)?
\item[\textbf{?R:}] \textit{How a fundamental competence is rule application in mathematics education?} Does respective technical support increase students' sense of safety in doing mathematics (e.g. when being doubtful whether steps in a proof are correct or not)?
\end{itemize}

\section{Specification Language --- a Means for Abstraction}\label{sec:spec-lang}
At the end of this section a more precise description  of ``applied mathematics'' will be clear: this is the part of mathematics, which is captured by formal specifications of problems associated by methods solving these problems and which also proves correctness of solutions. But first comes an example:
{\footnotesize\begin{tabbing}
1234\=123\=12\=12\=12,\=Postcond \=: \= $\forall \,A^\prime\, u^\prime 
\,v^\prime.\,$\=\kill
\>\texttt{01}\>Problem (Biegelinie, [Biegelinien])\\
\>\texttt{02}\>\> Specification:\\
\>\texttt{03}\>\>\> Model:\\
\>\texttt{04}\>\>\>\> Given  \>: ${\it Traegerlaenge}\;L, \;{\it Streckenlast}\;q_0$  \\
\>\texttt{05}\>\>\>\> Where  \>: $q_0 \; {\it ist\_integrierbar\_auf}\;[0,L]\;\land\;L>0$ \\
\>\texttt{06}\>\>\>\> Find   \>: ${\it Biegelinie}\;y$ \\
\>\texttt{07}\>\>\>\> Relate \>: ${\it Randbedingungen}\;[Q\,0=q_0\cdot L,\;M_b\,L=0, \;y\,0=0, \;\frac{d}{dx}y\,0=0]$\\
\>\texttt{08}\>\>\> References: \\
\>\texttt{09}\>\>\>\> Theory \>: Biegelinie \\
\>\texttt{10}\>\>\>x\> Problem \>: $[$"Baustatik", "Biegelinien"$]$ \\
\>\texttt{11}\>\>\>o\> Method \>: $[$"Integrieren", "KonstanteBestimmen"$]$ \\
\>\texttt{12}\>\> Solution:
\end{tabbing}}
\label{expl:spec-biegel}
\noindent
The line numbers \texttt{01..12} serve referencing in the subsequent paragraphs and do not belong to the specification. The example is taken from the description by the prototype mentioned~\cite{wn:proto-sys}. Further description will follow.

\subsection{Technical Facts and Prospects}\label{ssec:spec-facts}
\paragraph{A formal specification describes a problem most exactly:} The above example describes the problem of calculating a bending line (see~\footnote{An informal description of the problem is found at
\url{http://www.ist.tugraz.at/projects/isac/www/kbase/exp/exp_Statics_Biegel_Timischl_7-70.html}
within the prototype's example collection under \emph{Statics} at
\url{http://www.ist.tugraz.at/projects/isac/www/kbase/exp/index_exp.html}

}; the problem serves also as a running example in subsequent sections): The field \texttt{Given} above identifies the input required to solve the problem, \texttt{Where} is the pre-condition ensuring solvability, \texttt{Find} identifies the output expected. The field most characteristic for a problem is \texttt{Relate}, (an essential part of) the post-condition. It relates input, for instance $L,q_0$, and output, for instance $y(x) = \frac{q_0 \cdot L ^ 2}{4 \cdot EI} \cdot x ^ 2 + \frac{L \cdot q_0 }{6 \cdot EI} \cdot x ^ 3 + \frac{q_0}{24 \cdot EI} \cdot x ^ 4$ ($\it EI$ is a constant known in structural engineering). 

Formal specification of a mathematical problem is a prerequisite for reliably checking correctness: A solution like the function $y$ is correct, if substitution into the post-condition evaluates to true by an (automated) prover.

Associating specifications with one or more \texttt{Method}s solving the problem appears straight forward; the above example identifies a method in line \texttt{11}.

\medskip
More and more knowledge mechanised in proof assistants addresses problems from STEM faculties, for instance Networks, Security, Economics or Probability Theory in Isabelle's Archive of Formal Proofs~\footnote{\url{https://www.isa-afp.org/topics.html}}. But the kind of explicit specification addressed in the paper is still missing; so specifications need to be added to \emph{existing} knowledge --- assuming to have knowledge associated with specifications covering all undergraduate mathematics seems not unrealistic; see for instance~\footnote{\url{https://xenaproject.wordpress.com/2018/10/07/what-is-the-xena-project/}
}.

\paragraph{Trees of formal specifications allow automated problem refinement} as shown by the prototype already mentioned: \cite{richard:da} implements a tree for types of equations, which enables the system to determine the specification of the appropriate type of equation~\footnote{The prototype's sub-tree for types of equations can be found under \textit{equation} at
\url{http://www.ist.tugraz.at/projects/isac/www/kbase/pbl/index_pbl.html}}. Given a particular equation this is matched with $=$ without pre-conditions and then in a breadth-first search down the branches with tighter pre-conditions until the appropriate type of equation is found. This search resembles in a tracable refinement in comparison with computer algebra, which does the same by programme code (i.e. not tracable) in equation solvers.

\medskip
The large and active community in ``Formal Methods'' successfully works on ``systems of systems''~\cite{sys-of-sys-15} in close cooperation with industry; this might result in large collections of specifications in hierarchies of layers. And such collections of formal specifications might appear on workstations sooner or later in several engineering disciplines.

\begin{figure} [htb]
  \centering
  \includegraphics[width=115mm]{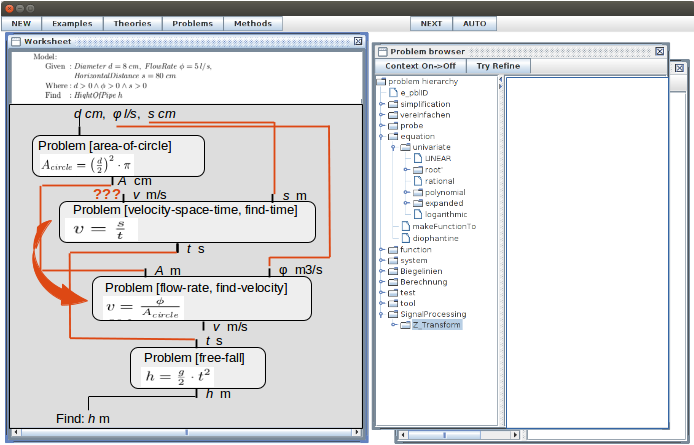} 
  \caption{Sequencing of specifications within ``divide and conquer''~.}
  \label{fig:sub-probl}
\end{figure}

\paragraph{Specifications become objects of interactive manipulation:} The already mentioned prototype investigated possibilities to regard specifications as objects for manipulation on the screen as shown in Fig.\ref{fig:sub-probl}.
This is a screenshot, where the panel on the left is faked, i.e. not yet existing; but the idea appears realisable, given specifications as described here together with the strong type systems in TP.

Fig.\ref{fig:sub-probl} shows a specification on top left and a panel underneath, where specifications can be placed as (black) boxes by ``drag and drop'' from the collection on the right. And then the task for students (or later on of engineers in more professional contexts) will be to arrange specifications such that respective inputs can be connected with respective outputs to a (directed and acyclic) graph of sub-problems to be tackled one by one in a subsequent calculation. The big red arrow indicates flipping two specifications which successfully completes the phase of problem solving in the example. A slide movie~\footnote{\url{http://www.ist.tugraz.at/projects/isac/publ/movie-sub-problems.pdf}}
tries to give an idea of the whole interactive process.

In the prototype user-guidance in modelling, i.e. in filling \texttt{Find} to \texttt{Relate}, is provided by minimal data (see \cite[p.~12]{wn:proto-sys} for more details) hidden behind the textual problem description. So, traditional example collections just need to add these data to make collections interactive.

\subsection{Foreseeable Effects and Questions}\label{ssec:spec-effects}
The previous section mentioned more prospects than facts; but the effects on education appear particularly interesting.

\paragraph{Separation of ``what'' from `` how'' in problem solving is supported:} An important principle in solving problems is, not only in engineering disciplines, to first clarify, \emph{what} the problem is, and then to consider \emph{how} to solve the problem --- this separation is now implemented by specifications separated from actual construction of solutions (discussed in \S\ref{sec:proof-lang} below).

Specifications are now objects the system can operate on (by automated refinement like with equations, by checking matches of output and input like in the slide-movie, etc) and the student can operate on: by selection from collections of specifications, by input to the fields \texttt{Given}..\texttt{Relate}, by arranging specifications in order to determine sub-problems, etc. Such operations promise to establish a new level of abstraction in mechanical support of problem solving. Such novelties raise questions to educators as well as to developers:

\begin{itemize}
\item[\textbf{?U:}] \textit{How guide students through the search space for specifications and methods} such that they don't get lost in large collections? Or: How introduce them to similar tasks on future workbenches for engineers?

\item[\textbf{?U:}] \textit{Which collections of specifications can be structured for automated refinement?} So far, only specifications for types of equations have been investigated~\cite{richard:da}. But apparently there are other areas of mathematics suitable for problem refinement and respective mechanisation --- a novel challenge for educators as well as for developers!
\item[\textbf{?U:}] \textit{Can switching between alternative solutions be supported smoothly?} Although alternative specifications and methods can pursue alternative ways towards solutions in different windows --- can interactive switching be supported in a way which is helpful for students?

\item[\textbf{?R:}] \textit{At which age can formal specification be introduced to education?} Early field tests with the prototype forced to develop a specific mode of interaction skipping specifications, because high-school teachers insisted, that specifications would disturb and overwhelm their pupils --- while specification (albeit not formal) is a mandatory item in high-school syllabi.

On the other side, academic educators in engineering disciplines give high priority to specifications; so there is a gap to be narrowed between high-school mathematics and academic courses.
\end{itemize}

\paragraph{Prerequisite knowledge can be determined automatically:}
Since an arbitrary selection of problems (in applied mathematics) to be solved or exercised are accompanied by respective formal specifications and since such a specification is associated with methods solving a respective problem class, and where methods in turn usually identify sub-problems (again accompanied by a specification and again associated with methods), all knowledge (down to elementary methods like simplification of fractions or integrating polynomials) required to solve the arbitrary selection can be determined automatically. In fact, the feature provides automated planning of learning paths~\cite{learn-path-09}. This raises questions.
\begin{itemize}
\item[\textbf{?R:}] \textit{How dynamically can curriculum planning be handled?} Since for any collection of problems in STEM studies all prerequisite knowledge can be calculated automatically (supposed, that respective knowledge is mechanised), this feature provides novel support for curriculum planning and raises further questions: What are the outcomes, if the feature is used to determine core knowledge for different kinds of schools? How far can go individual planning for counties, for schools, for classes?

\item[\textbf{?U:}] \textit{How improve independent learning by automated detection of prerequisite knowledge?} In case an individual student has chosen a problem (or a couple of problems) in a TP-based system, she or he can request presentation of all prerequisite knowledge --- how present the automatically created bunch of knowledge? How can that be made useful for students at which level?
\end{itemize}

\section{Proof Language --- Adapted to Mathematics for Engineers} \label{sec:proof-lang}
The level of students' mental maturity required to make teaching to prove reasonable is under debate. This paper bypasses this debate and pursues the idea to adapt TP technologies to mathematical traditions and not vice versa. The following is considered a calculation as close to traditions in mathematics for engineers as reasonable:
\label{expl:calc-biegel}
{\footnotesize\begin{tabbing}
123\=12\=12\=12\=12\=12\=12\=123\=123\=123\=123\=\kill
01\> Problem (Biegelinie, [Biegelinien]) \\
02\>\> Specification: \\
03\>\> Solution: \\
04\>\>\> Problem (Biegelinie, [vonBelastungZu, Biegelinien]) \\
05\>\>\> $[$\>$V\,x = c + -1\cdot q_0\cdot x, $\\
06\>\>\>    \>$M_b\,x = c_2 + c\cdot x + \frac{-1\cdot q_0}{2\cdot x ^ 2}, $\\
07\>\>\>    \>$\frac{d}{dx}y\,x =  c_3 + \frac{-1}{EI} \cdot (c_2\cdot x + \frac{c}{2\cdot x ^ 2} + \frac{-1\cdot q_0}{6\cdot x ^ 3}, $\\
08\>\>\>    \>$y\,x =  c_4 + c_3\cdot x +  \frac{-1}{EI} \cdot  (\frac{c_2}{2}\cdot x ^ 2 + \frac{c}{6}\cdot x ^ 3 + \frac{-1\cdot q_0}{24}\cdot x ^ 4)\;\;]$ \\
09\>\>\> Problem (Biegelinie, [setzeRandbedingungen, Biegelinien])\\
10\>\>\> $[$\>$L \cdot q_0 = c,\; 0 = \frac{2 \cdot c_2 + 2 \cdot L \cdot c + -1 \dot L ^ 2 \cdot q_0}{2},\; 0 = c_4,\; 0 = c_3\;\;]$\\
11\>\>\> solveSystem $(L \cdot q_0 = c,\; 0 = \frac{2 \cdot c_2 + 2 \cdot L \cdot c + -1 \cdot L ^ 2 \cdot q_0}{2},\; 0 = c_4,\; 0 = c_3\;\;],\; [c, c_2, c_3, c_4])$\\
12\>\>\> $[$\>$c = L \cdot q_0 ,\; c_2 = \frac{-1 \cdot L ^ 2 \cdot q_0}{2},\; c_3 = 0,\; c_4 = 0]$\\
13  \` Take $y\,x =  c_4 + c_3\cdot x +  \frac{-1}{EI} \cdot  (\frac{c_2}{2}\cdot x ^ 2 + \frac{c}{6}\cdot x ^ 3 + \frac{-1\cdot q_0}{24}\cdot x ^ 4)$ \\

14\>\>\> $y\,x = c_4 + c_3 \cdot x + \frac{-1}{EI} \cdot (\frac{c_2}{2} \cdot x ^ 2 + \frac{c}{6} \cdot x ^ 3 + \frac{-1 \cdot q_0}{24} \cdot x ^ 4)$\\
15  \` Substitute $[c,c_2,c_3,c_4]$ \\

16\>\>\> $y\,x = 0 + 0 \cdot x + \frac{-1}{EI} \cdot (\frac{\frac{-1 \cdot L ^ 2 \cdot q_0}{2}}{2} \cdot x ^ 2 + \frac{L \cdot q_0}{6} \cdot x ^ 3 + \frac{-1 \cdot q_0}{24} \cdot x ^ 4)$\label{exp-biegel-Substitute}\\
17  \` Rewrite\_Set\_Inst $([({\it bdv},x)], {\it make\_ratpoly\_in})$ \\

18\>\>\> $y\;x = \frac{q_0 \cdot L ^ 2}{4 \cdot EI} \cdot x ^ 2 + \frac{L \cdot q_0 }{6 \cdot EI} \cdot x ^ 3 + \frac{q_0}{24 \cdot EI} \cdot x ^ 4$\\
19\> $y\;x = \frac{q_0 \cdot L ^ 2}{4 \cdot EI} \cdot x ^ 2 + \frac{L \cdot q_0 }{6 \cdot EI} \cdot x ^ 3 + \frac{q_0}{24 \cdot EI} \cdot x ^ 4$
\end{tabbing}}
\noindent
The line numbers \texttt{01..19} serve referencing in this paper and do not belong to the example. The steps constructing a \texttt{Solution} are on the left; indentation indicates that a term takes more than one line (e.g. in \texttt{05..08}) or indicates that the subsequent lines with deeper indentation could be collapsed (like sub-directories in a file browser): a click on line \texttt{03} would collapse lines \texttt{04..18} as shown on p.\pageref{expl:spec-biegel} and a click on line \texttt{02} above would expand the \texttt{Specification} shown on p.\pageref{expl:spec-biegel} and there on line \texttt{12}. Further description follows subsequently.

\subsection{Technical Facts and Prospects}\label{ssec:proof-facts}
\paragraph{Forward proof can be adapted to calculations in applied mathematics} as shown above: the \texttt{Solution} of the problem to calculate the bending line of a beam under load is considered as close to what would be written on a blackboard traditionally --- with the difference already mentioned, that there are no long arrows to asides on the blackboard, rather, software allows to exhibit and hide details in a structured manner. So let's skip the first lines of the \texttt{Solution}, leave it to the subsequent paragraph on modularisation and have a look at lines \texttt{13..18}. On the left there are the steps of calculation leading to the function $y\;x = \frac{q_0 \cdot L ^ 2}{4 \cdot EI} \cdot x ^ 2 + \frac{L \cdot q_0 }{6 \cdot EI} \cdot x ^ 3 + \frac{q_0}{24 \cdot EI} \cdot x ^ 4$ solving the problem. On the right there are tactics, which justify the respective step. Line \texttt{13} tells which term to continue with, \texttt{15} tells the variables to be substituted from the context ($c,c_2,c_3,c_4$ have been calculated in \texttt{12} before) and \texttt{17} is still somewhat cryptic (with apologies for a prototype's fragmentariness) for an elementary task, which is simplification of the term in \texttt{16} by use of the rule-set \texttt{make\_ratpoly\_in} for rational polynomials with  $x$ as bound variable, \texttt{bdv}. A specific click on \texttt{16} would show the intermediate steps from $y\,x = 0 + 0 \cdot x + \dots$ to the simplified result.

This format is what has been introduced as ``structured derivation'' \cite{back-SD-2010} long time ago. The format is based on natural deduction~\cite{Gentzen:nat-deduct}, the logical foundation of most proof assistants. Natural deduction is the meta-language for the proof language Isabelle/Isar~\cite{wenzel:isar}, too, where the latter is human readable and the former is program code~\footnote{Program code implementing the laws of natural deduction in Isabelle is called the ``trusted kernel'' and kept together in one file~\url{https://isabelle.in.tum.de/repos/isabelle/file/91162dd89571/src/Pure/thm.ML}}.

The logical foundations ensure the calculation's ``correctness by construction'' step by step; and correctness of the result can be (automatically) proven by substitution into the post-condition on p.\pageref{expl:spec-biegel} line \texttt{07}. The format resembles ``forward proof'' in proof assistants with the difference, that the system does not insist to make the result already explicit at the beginning (as, for instance, done by the \texttt{lemma} in Fig.\ref{fig:isa-transparent}).

\emph{Interactive} construction of a calculation like above will be discussed in \S\ref{sec:next-step}.

\paragraph{Calculations can be modularised} according to the principle ``divide and conquer'' in problem solving. The above example-\texttt{Solution} contains three sub-problems in line \texttt{04}, \texttt{09} and \texttt{11} respectively, all immediately followed by the respective result. The latter sub-problem is a call for an equation solver similar to computer algebra.  \S\ref{ssec:spec-facts} showed (remember Fig.\ref{fig:sub-probl}), how the sequence of sub-problems can be interactively determined, how data can be input to the model's fields \texttt{Given}..\texttt{Relate} and how methods solving the problems can be selected from libraries. 

Entering a method marks the transition from modelling / specifying a problem to step-wise constructing a solution of the problem; so we have nested recursion on the phases of modelling, specifying and solving implemented within one (prototype-)system. For specification students use the  \texttt{References}' fields \texttt{Theory}, \texttt{Problem} and \texttt{Method} on p.\pageref{expl:spec-biegel}. The \texttt{x} in line \texttt{10} indicates, that the \texttt{Problem}'s \texttt{Model} is shown (and not the \texttt{Method}'s guard, line \texttt{11}).

\paragraph{The power of TP is accompanied by adaptive user-guidance.} \S\ref{ssec:term-facts}..\S\ref{ssec:proof-facts} gave a presentation of technology and \S\ref{ssec:next-facts} will show, how that technology allows for flexible students' interaction far beyond the state of present educational software in mathematics. Managing flexible interaction calls for a dialogue module; the design of such a module attracted much efforts in prototype development, e.g.~\cite{kremp.np:assess,AK04:thesis}. Detailed design and implementation of rules for the dialogues~\cite{mkienl-bakk} concerns a major development effort in collaboration between computer mathematics and educational sciences. Prospects assume a complexity in the assembly of these rules such that another separate language layer appears on the horizon, the layer of dialogue rules creating a separate development task for ``dialogue authors''.

\subsection{Foreseeable Effects and Questions}\label{ssec:proof-effects}
What have been step-wise calculations on a blackboard or in a notebook on paper, that has become interactive step-wise constructions supervised by a mechanised system. This is a fundamental change with remarkable effects on education.

\paragraph{Students can trust that their problem solutions are correct:} 
The use of computer algebra led to situations overwhelming students: using a mechanical system even increases difficulties to decide, what is correct or not. In mathematics the decision for being correct or not is essential, so this paper assumes that deciding correctness within formal mathematics is \emph{not} in responsibility of students, but of TP-based systems. We see the following questions:
\begin{itemize}
\item[\textbf{?R:}] \textit{Do systems' checks for correctness foster students' self-assuredness or students' credulity?} Both effects can be imagined, so: How design software of the new generation such that students' self-assuredness is developed and not credulity? Or asked more generally: How can mechanised mathematics systems decrease impressions of magic in mathematics?

\item[\textbf{?U:}] \textit{Which possibilities are there for immediately handling errors in user input?} This is an interesting question with respect to technology which is designed for ``correctness by construction''. One possibility is to react on ``error patterns'' in simplification~\cite{gdaroczy-EP-13} which has already be prototyped.

\item[\textbf{?U:}] \textit{Can input of wrong steps be handled until construction of steps reaches contradictions?} Isabelle/Isar~\cite{wenzel:isar} provides a \texttt{sorry} in order to post-pone sub-proofs and first clarify other parts in formal deduction; such a feature appears desirable within mathematics for engineers --- and then the above question is asking for how far such a feature can go in structured derivations.
\end{itemize}

\paragraph{Students can inspect justifications down to elementary steps:} 
Proof assistants got their name, because they serve interactive proof construction, where a mathematician constructing a proof has justifications at different levels of abstraction in mind. If he or she has implemented the justifications, these can be inspected, of course. Most of the rules are applied by use of matching and rewriting as introduced in \S\ref{ssec:spec-facts}. 
If Isabelle/Isar's feature of transparency is maintained during adaption to mathematics for engineers, the following questions arise:

\begin{itemize}
\item[\textbf{?U:}] \textit{How present simplifications comprising many dozens of rewrites?} Isabelle's simplifier got tracing features~\cite{hupel-simp-trace}, triggered by users' requests for more transparency; however, the implementation addresses expert users and is not immediately useful for novices. Thus there are trials in prototyping~\cite{MG:thesis} to be reviewed.

\item[\textbf{?U:}] \textit{How present answers to questions about basic reasoning steps in natural deduction?} As mentioned in \S\ref{ssec:proof-facts} the laws of natural deduction are written in a program language assumed not readable for students in mathematics education. So this question is an issue, but not urgent, since experiences show that structured derivations appear evident to students without going down to natural deduction.
\end{itemize}

\paragraph{All phases of problem solving are covered by \emph{one} software:} 
This feature allows for independent learning: given a problem description by text or figures, the first phase of problem solving is modelling, i.e. translation into terms to be input to the fields \texttt{Given}..\texttt{Find}; the second phase is specifying, i.e. relating the model to an appropriate specification in the knowledge base and selecting an appropriate method; the next phase is to construct a solution step by step and the final phase is an automated check, if the post-condition is fulfilled, i.e. that the solution is correct. These questions are raised:

\begin{itemize}
\item[\textbf{?R:}] \textit{What happens when students tackle completely new problems?} An extreme view could be: Given a systems of the new generation and a problem implemented in such a system, no knowledge need to be presupposed from a student, because she or he can ask the system step by step for all knowledge required to solve the problem. So, in which situations and in which points is this extreme view unrealistic?

\item[\textbf{?U:}] \textit{How adapt user-guidance between different phases of problem solving?} With respect to technology, user-guidance is different in different phases involving different language layers. So, how make user-guidance smooth and uniform over all phases, adapting to various users from novices to experts, from students experimenting with ideas to those who only want to exercise, etc?
\end{itemize}

\paragraph{Software can connect introductory and advanced courses:} 
There is a principal and frequently observed problem in mathematics education: introductory courses (or high-school mathematics in general) provide prerequisites for advanced courses, but in higher semesters of STEM studies much has been forgotten from introductory courses. One reason for the phenomenon seems obvious: Introductory courses lack convincing application of what is being taught, with at least two good reasons: (1)  many courses collect students from various engineering disciplines, (2) the students lack advanced knowledge to understand convincing applications. As a consequence relevance of introductory mathematics for specific engineering studies is not acknowledged by students and thus these are not motivated to learn and keep in mind.

Now flexibility of TP-based software allows to present \emph{one and the same} collection of examples to students of both, of introductory courses and of advanced courses, if provided adaptive user-guidance (see \S\ref{ssec:proof-facts}). In the former the basic method just introduced is in the focus of the dialog --- while the other parts of the calculations are open for investigation; in the latter specification and surveys on sub-problems are in the focus --- while details of basic methods probably forgotten can be re-done.

With such software in the back, a lecturer can say in introductory courses: For the method (e.g. integration) there is no time to demonstrate serious application, but you will see lots of applications, if you investigate the surrounding parts in the calculations! And in the advanced courses can be said: Don't worry, if you have forgotten some basic methods (e.g. integration), we don't need time to repeat these --- just go into the details of the obligatory exercises and review what you have forgotten!

\begin{itemize}
\item[\textbf{?U:}] \textit{How change user-guidance from introductory courses to advanced courses?} Apparently, user-guidance must be different, if a part of a calculation concerning a specific basic method (e.g. integration) is picked out for exercising (in introductory courses) --- or if many sub-problems are handled as black-boxes (e.g. integration) and surveys on dependencies between sub-problems are relevant (in advanced courses).

\item[\textbf{?U:}] \textit{How decide between accepting black-boxes and enforcing active construction from students?} Students can calculate examples with various intentions, for instance, switching between detail and survey as described above for different courses. How design user-guidance such that it adapts also to intentions of lecturers appropriately?
\end{itemize}

\section{Programming Language --- for Easy Authoring}\label{sec:prog-lang}
The layer of programming languages is not visible for students, but for authors of problem collections, who want to contribute specifications of problems and methods solving the latter. Since invisible for students, we skip effects for education and address technical facts only. The latter are also required as a prerequisite for \S\ref{sec:next-step}, which will exhibit educational relevance again.

\begin{figure} [htb]
  \centering
  \includegraphics[width=145mm]{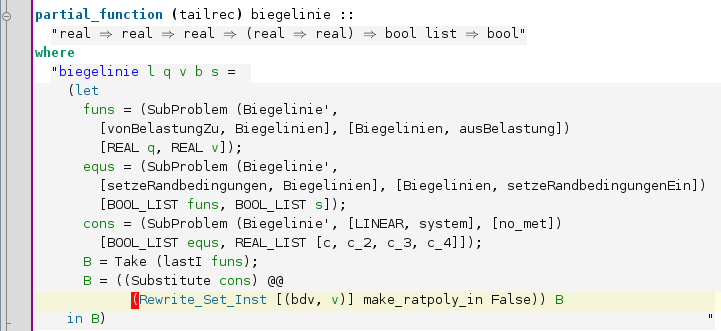}
  \caption{New appearance of the program introduced in \cite[p.~92]{wn:proto-sys}.}
  \label{fig:fun-pack-biegelinie}
\end{figure}

Computer algebra has been equipped with programming languages for a long time. For TP an integration of programming languages is rather new, see for instance Isabelle~\cite{krauss}. Challenges for such programming are described in~\cite{plmms10}. These have not changed so far, except that prototyping already started to shift programmes into Isabelle~\cite{wn:lucin-thedu18}. Such a program is shown in Fig.\ref{fig:fun-pack-biegelinie}.
The program \texttt{biegelinie} takes five arguments \texttt{l, q, v, b, s}; the types of these arguments are made explicit in the signature above. In the program body first come the three \texttt{SubProblem}s also found in the calculation on p.\pageref{expl:calc-biegel}. Then the tactic \texttt{Take} appears in the program like on p.\pageref{expl:calc-biegel} (in line \texttt{13} at the right margin). Subsequently two functions are concatenated by \texttt{@@} which both concern tactics also shown (by user request) in the calculation.

The following technical facts and prospects concern this kind of programs.

\paragraph{Algorithms are programmed similar to computer algebra}
with essential differences in the logical background. First, the programs are purely functional, there are no side effects like input or output --- which enhances provability. A second difference is, that Isabelle's logic requires proofs of termination for all programs. However, prototyping showed that generality of most programs does not allow simple termination proofs; so this issue is bypassed by implementing \texttt{partial\_function}s without proof. A third difference, the most essential one, is that interpretation each tactic updates a logical context (for instance encountering a term $\frac{\dots}{x-1}$ and recording $x\not=1$).

Many issues are common to TP and to computer algebra, for instance, multi-valuedness of algebra results or real numbers and approximation~\cite{plmms10}, resolving these will take time but appears not urgent for educational purposes.

Another important difference to programs in computer algebra is, that function calls are replaced by interactive modelling and specifying sub-problems by the user. Thus a \texttt{SubProblem} takes for arguments: a theory containing the key to respective knowledge (e.g. \texttt{Biegelinie} in Fig.\ref{fig:fun-pack-biegelinie}), a specification addressed by a key (e.g. \texttt{[vonBelastungZu, Biegelinien]}), a method addressed by another key (e.g. \texttt{[Biegelinien, ausBelastung]}) and a list of formal arguments (e.g. \texttt{[q, v]}). The program describing how to solve the respective sub-problem has a guard with structure analogous to the \texttt{Model} on p.\pageref{expl:spec-biegel}. The guard prevents the program from being executed in case the pre-conditions are not fulfilled. Where interaction comes from in a functional program, this is introduced in the next paragraph and discussed in more detail in \S\ref{ssec:next-facts} below.

\paragraph{Programming mathematics is separated from programming interaction.} The program in Fig.\ref{fig:fun-pack-biegelinie} is purely functional, as already mentioned, so: Where does interaction, output and students' input come from, which is indispensable in an educational system? The answer will be given in the next \S\ref{sec:next-step}. Here we state an important separation of concerns: mathematical algorithms are implemented in the program language under consideration such that user interaction needs not to be considered; mathematicians can focus mathematical algorithms.

Interaction is concern of educators and psychologists. These are freed from mathematical details and can focus development of support for educational processes. As mentioned in \S\ref{ssec:proof-facts}, a separate language layer will emerge. Elements of this language need to be defined in a user model, which needs to be designed from scratch, since the features of TP-based systems are beyond what could be imagined so far. 

\paragraph{Programs are general and re-usable.} 
The program in Fig.\ref{fig:fun-pack-biegelinie} looks very specific, but it solves all problems found on two pages in a textbook~\footnote{\url{http://www.ist.tugraz.at/projects/isac/www/kbase/exp/exp_Statics_Biegel_Timischl.html} as a part of the prototypes problem collection \url{http://www.ist.tugraz.at/projects/isac/www/kbase/exp/index_exp.html}.}. Nevertheless, much effort is required to implement all the knowledge required to make a TP-based system ready for education: programs, problem specifications and other definitions. However, once the basic methods are programmed and respective specifications are formalised, authoring of new problems should be efficient like in computer algebra, where libraries are at programmers' disposal.

\section{``Next-Step-Guidance'' --- a Novel Contribution}\label{sec:next-step}
The novel contribution to TP-technology for mathematics education is combining computation and deduction in a way, which is described in \cite{wn:lucas-interp-12} with respect to logical foundations and which has been thoroughly investigated in the prototype. It has been called ``Lucas-Interpretation'' in honour of a mentor of the design and prototyping process.

\subsection{Technical Facts and Prospects}\label{ssec:next-facts}
Details of technology from users' perspective are described in \cite{thedu16:lucin-user-view}, with respect to this paper has to be stated: Lucas-Interpretation involves all language layers mentioned above: Terms from \S\ref{sec:term-lang} are handled as elements of calculations, specifications from \S\ref{sec:spec-lang} are manipulated interactively in case sub-problems are encountered in a calculation, step-wise construction of a calculation follows logic foundations of proofs from \S\ref{sec:proof-lang} and the interpreted programs belong to the language from \S\ref{sec:prog-lang}.

\paragraph{Lucas-Interpretation can determine a next step in problem solving.}
Interpretation operates on programs described above in \S\ref{sec:prog-lang}. That means, this kind of next-step-guidance is provided during step-wise construction of solutions for a given problem (while next-step-guidance during modelling and specifying relies on other technologies).

The ``\emph{determine} a next step'' concerns the \emph{calculational} part of the interpreter, alien to principles of proof assistants (where general programs constructing proofs are not expected): A program just determines a next step within an algorithm constructing solutions for a specific class of engineering problems. The ``\emph{can} determine'' indicates, that the request for a next step can come either from a student \emph{or elsewise} from the dialog module introduced in \S\ref{ssec:proof-facts}: If a student gets stuck, she or he can request a next step. But if the user model says, for instance, that the student is in exercise mode and that he or she already caused series of similar requests, then the dialog module might decide to reject the request --- and to request more activity from the student in turn: input a term or at least select a rule from a list to apply, etc.

\paragraph{Lucas-Interpretation allows to prove correctness of students' input} during step-wise construction of solutions for a given problem. The ``allows to \emph{prove}'' concerns the \emph{deductive} part of the interpreter, which is harder than the computational part: 

In a calculation a step input by a student is either a term (or a part of a term, see \S\ref{ssec:proof-facts}) or a tactic. The latter case is the simpler one: the interpreter has to search for a respective tactic in the program code. Difficulties here only arise with partial equivalence of tactics or with algebra rules not contained in a rule set. Experiences with the prototype indicate, that such difficulties can be overcome by careful programming.

In the case a student inputs a term, the interpreter is challenged to find an automated prover able to derive the term from the current context in the calculation. The present prototype only uses a simplifier and, nevertheless, shows good results. The generated derivation is inserted into the calculation as intermediate steps. These are one level deeper in the calculation, thus hidden and only displayed on request.


\subsection{Foreseeable Effects and Questions}\label{ssec:next-effects}
Here effects heavily depend on the dialogue module introduced in \S\ref{ssec:proof-facts}. The design of this module anticipates two general requirements for user modelling in educational systems: (1) Maintain dialogues in a balanced way, neither overwhelming nor boring for an individual student and (2) adapt to various requirements for learning raised by students as well as by course instructors. Below these points are addressed in sequence.

\paragraph{Equitable dialogues between student and system.} The technology of Lucas-Interpretation opens up much of the power and potential of TP technology. The power appears sufficient to strive for ``dialogues on an equal base'' between student and system: both of them can solve problems and have respective knowledge in mathematics, the system shall be polite enough to listen to the student until something is wrong. Then the system will open up all knowledge requested by the student; but also the system can ask the student, for instance to apply a rule or to select a method appropriate to solve a specified problem, etc. And the system shall monitor a balance in students' interaction: as active as possible, but not overwhelming and never boring! Detailed design of a dialog module and a respective user model raises questions like

\begin{itemize}
\item[\textbf{?U:}] \textit{How detect when a student gets bored?}
\item[\textbf{?U:}] \textit{How determine when a student gets overwhelmed?}
\item[\textbf{?U:}] \textit{How distinguish silent thinking from other kinds of breaks?}
\item[\textbf{?:U}] \textit{\dots}
\end{itemize}

\paragraph{Dialogues adapt \dots}$\;$\\
\dots several situations already addressed in previous sections. Since now all the language layers and respective technical facts as well as foreseeable effects have been presented, we give a brief survey of the many questions arising from user-modelling in a more structured way.
\begin{itemize}
\item \textbf{\dots adapt to levels of knowledge} which a student already has worked through in solving problems more or less successfully (which has been recorded in a user model).
  \begin{itemize}
  \item[\textbf{?U:}] \textit{How determine success in solving problems} given in certain collections of problems?
  \item[\textbf{?U:}] \textit{How decide that a student is able to solve a certain kind of problem?} To what extent?
  \item[\textbf{?U:}] \textit{How simplify pre-configuration of black-boxes by course instructors?} Certain kinds of problems, for instance calculating derivatives, can be declared as ``black-box'' such that the system skips interactive solving, just presents the respective result and continues silently until the final result.
  \item[\textbf{?U:}] \textit{Under what circumstances re-open black-boxes} and enforce interactivity, e.g. for repetition?
  \item[\textbf{?R:}] \textit{How appropriate are existing user models} for TP-based educational systems?
  \end{itemize}

\item \textbf{\dots adapt to students' purpose to use the system,}  which can be: exercise a certain collection of problems, exercise a certain method in various applications, explore a collection of problems which are new to students, find applications of a certain method in various areas, etc.
  \begin{itemize}
  \item[\textbf{?U:}] \textit{How determine a student's purpose of usage?} By explicit input and/or by inference from ongoing interaction? The latter would meet the frequent experience, that students not exactly know what they want or need in specific situations.
  \item[\textbf{?U:}] \textit{What are the system's features appropriate to support a certain purpose?}
  \item[\textbf{?R:}] \textit{What are the discriminants to identify different purposes of usage?}
  \end{itemize}

\item \textbf{\dots adapt to a focus of instruction} which needs not be analogous to students' purposes mentioned above. The side of instruction raises other questions:
  \begin{itemize}
  \item[\textbf{?U:}] \textit{How graceful handle students overriding a focus} set by instructors?
  \item[\textbf{?U:}] \textit{How design course management} such that a focus can be set easily?
  \item[\textbf{?R:}] \textit{Which instructional modes can be identified} for capturing above kinds of focus?
  \item[\textbf{?:}] \dots
  \end{itemize}

\item \textbf{\dots adapt to individual preferences of perception} like visual, auditory, haptic, etc.
  \begin{itemize}
  \item[\textbf{?U:}] \textit{\dots}
  \end{itemize}
\item \textbf{\dots }
\end{itemize}
\dots lots of challenging questions arising from power, flexibility and generality of TP-based systems: since these provide features far beyond the state of the art, efforts also on the educational side seem justified, in order to bring strengths of technology to bear in educational practice --- as said already, a challenge for collaboration between computer mathematics and educational sciences.


\section{Summary and Conclusions}\label{sec:summary}
This study collected specific technologies from computer theorem proving (TP) and from TP-based prototypes --- specific technical features, which appear relevant for developing educational software for mathematics at the state of the art. A major part of the features has been implemented and demonstrated during prototyping (in the ISAC-project and others, references are given), other features are in a concrete stage of planning.

The collection and statement of technologies and features serves the purpose to inform educators and to invite for collaborative development, see the conclusions below. The paper prepares for such collaboration by questions structured alongside four layers of formal language, combined to a ``complete, transparent and interactive model of mathematics'': (1) the language of mathematical terms as the basis for the other layers --- with types implementing abstract mathematical structures and with definitions such that they are open for students' investigation by mouse-click in \S\ref{ssec:term-effects}, featuring ``transparent'' models, which are also ``complete'' with respect to deduction. (2) the language of formal specifications --- determining problems in mathematics and in engineering, and a means for supporting ``divide and conquer'' in problem solving in \S\ref{ssec:spec-effects} featuring ``complete'' models covering \emph{all} phases in problem solving beginning with translation from problem statements into terms. (3) the language of proofs adopted from ``calculational proofs'' in proof assistants --- presented in a format which enhances traditional notation found on blackboards in academic lectures as well as in schools in \S\ref{ssec:proof-effects} featuring ``interactive'' as well as ``transparent'' models with respect to justification of steps towards problem solutions. (4) the language for programming algorithms for solving classes of problems --- an offer to lecturers and teachers to implement their preferred examples and exercises in \S\ref{sec:prog-lang}.

An additional section presents Lucas-Interpretation for implementing ``next-step-guidance'' --- a feature indispensable for educational systems, which should be able to propose next steps to students who got stuck in interactive problem solving in \S\ref{ssec:next-effects}. This feature is considered the final cornerstone for ``complete, transparent and interactive models of mathematics'' for education in STEM studies.

\paragraph{Conclusion on cooperation computer mathematics --- educational sciences:} The author hopes that the collection of technical features in this paper is illustrative enough to motivate educators to some of the following steps:
\begin{enumerate}
\item Review existing prototypes for the features described in this paper.
\item\label{Refine} Refine the questions posed in this paper with respect to appropriate frameworks in educational theories. The questions are distinguished into two strands: one is considered to address genuine research in educational sciences, another is biased towards software usability and to be solved in close collaboration with computer mathematicians. Some questions appear rather general, nevertheless raised again by challenges due the power of TP.
\item\label{Collaborate} Collaborate with computer mathematicians in developing prototypes to a stage appropriate for experimental use in classroom. Small field tests have been mentioned as successful in \S\ref{sec:intro}, but the tests also showed that usability in classroom requires considerable efforts to develop prototypes further to such a stage. 
\item\label{Design} Design extensive field tests inside and outside of classes at various levels. The various levels then should be related to the refined questions from Pt.\ref{Refine}. The challenge with such research questions seems, that the technical features described appear to support combinations of various parallel tasks in doing mathematics.
\item Proceed in cycles through Pt.\ref{Collaborate}..Pt.\ref{Design} establishing a new generation of TP-based software for mathematics within STEM education.
\end{enumerate}


\medskip
A possible approach from the side of educational sciences could be to imagine ``systems that explain themselves'' announced in \cite{wn:proto-sys}, somewhat naively from an educational point of view. From the side of technology seems clear: \emph{If possible at all, self-explanatory systems are possible in mathematics: the matter is abstract and less related to physical matter than in all other sciences. Thus relations remain within abstractions themselves, which can completely be captured by symbols implemented in computer software and operations on these symbols with automatically generated feedback.}

Naivety needs to be overcome on both sides: By computer mathematicians expressed as ``We can tell what the upcoming software generation is good for in education!'' and in educational sciences expressed as ``We have all we need'' or as ``Technology is irrelevant for proper mathematics education'' or as ``Mathematical software is part of educational problems rather than part of solutions''.


\bibliographystyle{eptcs}
\bibliography{references}

\end{document}